\documentclass[a4paper, 12pt]{article}

\usepackage{amsmath,amssymb}
\usepackage[numbers]{natbib}
\usepackage{rotating}
\usepackage{graphicx}

\newcommand{\R}{\mathbb{R}}
\newcommand{\N}{\mathbb{N}}

\begin{document}

\title{Existence Of Entire Solutions For Singular \\ Quasilinear Convective Elliptic Systems}

\author{
	\bf Umberto Guarnotta\\
	\small{Dipartimento di Matematica e Informatica, Universit\`a di Catania,}\\
	\small{Viale A. Doria 6, 95125 Catania, Italy}\\
	\small{\it E-mail: umberto.guarnotta@phd.unict.it}
}
\date{}
\maketitle

\begin{abstract}
The existence of entire solutions to quasilinear elliptic systems exhibiting both singular and convective reaction terms is discussed. An auxiliary problem, obtained by `freezing' the convection terms and `shifting' the singular ones, is first solved. Then, a priori estimates, fixed point arguments, nonlinear regularity, compactness results concerning the gradient terms, besides a regularization-localization procedure, yield the existence of a weak solution to the problem.
\end{abstract}

\section{Introduction}

Let us consider the quasilinear elliptic system
\begin{equation}
\label{problem}
\left\{
\begin{array}{lll}
-\Delta_p u = f(x,u,v,\nabla u,\nabla v) \quad &\mbox{in} \;\; &\R^N, \\
-\Delta_q v = g(x,u,v,\nabla u,\nabla v) \quad &\mbox{in} \;\; &\R^N, \\
u,v>0 \quad &\mbox{in} \;\; &\R^N,
%u(x), \, v(x) \to 0 \quad &\mbox{when} \; \; &|x| \to +\infty,
\end{array}
\right.
\end{equation}
with $ N \geq 3 $ and $ 1 < p,q < N $. We suppose $ f,g: \R^N \times (0,+\infty)^2 \times \R^{2N} \to (0,+\infty) $ to be Carathéodory functions satisfying the following growth conditions:
\begin{eqnarray*}
&m_1 a_1(x) s_1^{-\alpha_1} s_2^{\beta_1} \leq f(x,s_1,s_2,{\bf t}_1,{\bf t}_2) \leq M_1 a_1(x) (s_1^{-\alpha_1} s_2^{\beta_1} + |{\bf t}_1|^{\gamma_1} + |{\bf t}_2|^{\delta_1}), \\
&m_2 a_2(x) s_1^{\alpha_2} s_2^{-\beta_2} \leq g(x,s_1,s_2,{\bf t}_1,{\bf t}_2) \leq M_2 a_2(x) (s_1^{\alpha_2} s_2^{-\beta_2} + |{\bf t}_1|^{\gamma_2} + |{\bf t}_2|^{\delta_2}),
\end{eqnarray*}
for every $ (x,s_1,s_2,{\bf t}_1,{\bf t}_2) \in \R^N \times (0,+\infty)^2 \times \R^{2N} $, where $ m_i, M_i > 0 $, $ 0 \leq \alpha_1, \beta_2 < 1 $, $ 0 \leq \beta_1, \delta_i < q-1 $, $ 0 \leq \alpha_2, \gamma_i < p-1 $, and the weights $ a_i $ satisfy $ {\rm ess \, inf}_{B_\rho} \, a_i > 0 $ for any $ \rho > 0 $. \\
In order to get Sobolev and $ L^\infty $ estimates\footnote{Sobolev estimates furnish a control on the norms $ \|\nabla u\|_{L^p(\R^N)} $ and $ \|\nabla v\|_{L^q(\R^N)} $, while the $ L^\infty $ estimates control $ \|u\|_{L^\infty(\R^N)} $ and $ \|v\|_{L^\infty(\R^N)} $. From now on, otherwise explicitly stated, estimates are assumed to be global, that is, uniform on the whole $ \R^N $.} on the solutions of (\ref{problem}), representing the standard way to gain compactness, we are compelled to make further assumptions on both the growth exponents and the weights, namely,
\begin{equation}
\label{mixedcond}
\max\{\beta_1,\delta_1\} \max\{\alpha_2,\gamma_2\} < (p-1-\gamma_1)(q-1-\delta_2),
\end{equation}
to ensure a Sobolev bound, and
\begin{equation}
\label{weightsum}
a_i \in L^1(\R^N) \cap L^{\zeta_i}(\R^N), \quad \zeta_i \in (N,+\infty],
\end{equation}
with
\begin{eqnarray*}
\frac{1}{\zeta_1} < 1-\frac{p}{p^*} -\theta_1,& &\theta_1 := \max \left\{ \frac{\beta_1}{q^*},\frac{\gamma_1}{p},\frac{\delta_1}{q} \right\} < 1 - \frac{p}{p^*}, \\
\frac{1}{\zeta_2} < 1-\frac{q}{q^*} -\theta_2,& &\theta_2 := \max \left\{ \frac{\alpha_2}{p^*},\frac{\gamma_2}{p},\frac{\delta_2}{q} \right\} < 1 - \frac{q}{q^*},
\end{eqnarray*}
to obtain a $ L^\infty $ bound; here, $ p^* := \frac{Np}{N-p} $. It is worth noticing that (\ref{mixedcond}) is a mixed condition, that is, it takes into account the data of both equations in (\ref{problem}); no other mixed conditions occur. We also point out that, in the non-singular case $ \alpha_1 = \beta_2 = 0 $, the required properties on $ \zeta_1 $ (and analogously for $ \zeta_2 $) appearing in (\ref{weightsum}) guarantee that the right-hand side of (\ref{problem}) belongs to $L^r(\R^N) $, $ r > \frac{N}{p} $, which is the minimum requirement (among Lebesgue spaces) on the right-hand side of the $p$-Poisson equation $ -\Delta_p w = h(x) $ to get $ w \in L^\infty(\R^N) $, making condition (\ref{weightsum}), in a certain sense, natural. \\
The prototype of (\ref{problem}), obtained by setting $ \gamma_i = \delta_i = 0 $ and $ m_i = M_i $, has a \textit{cooperative structure}, i.e., $ f $ is increasing in $ v $ and $ g $ is increasing in $ u $; however, we require no monotonicity assumptions on $ f,g $. The Dirichlet version of (\ref{problem}) in bounded domains has been investigated in \cite{CLM}, while \cite{MMM} deals with (\ref{problem}) for $ \alpha_2 = \beta_1 = 0 $ and without convection terms (i.e., terms depending on the gradient of solutions). \\
The present investigation follows the direction of the recent papers \cite{LMZ,GMMot,GM}, regarding singular convective problems in bounded domains, with different boundary conditions. The literature about singular problems and convective ones is very wide: here we limit ourselves to addressing the reader to the monograph \cite{GR}, concerning mainly semilinear equations (vide Sections 4.7 and 9.4); see also the short presentation \cite{G} and the references therein. Incidentally, the results discussed here were obtained, under stronger hypotheses, in \cite{GMMou}; the extensions presented here are contained in a work in progress.

\section{The technique}

\bigskip

\noindent
\textbf{Step 0: The functional setting.}

\medskip

We look for solutions to (\ref{problem}) in a suitable functional setting, which is given by Beppo Levi spaces, also called homogeneous Sobolev spaces. Before introducing them, we briefly explain the main reason that led us to consider these spaces (following \cite[p. 80]{Galdi}). Let us consider the Dirichlet problem
\begin{equation}
\label{exterior}
\left\{
\begin{array}{lll}
\Delta u = 0 \quad &\mbox{in} \;\; &\R^3 \setminus \overline{B}_1, \\
u=1 \quad &\mbox{on} \;\; &\partial B_1, \\
u(x) \to 0 \quad &\mbox{when} \; \; &|x| \to \infty,
\end{array}
\right.
\end{equation}
being $ B_1 $ the ball of unitary radius centered at the origin. Setting $ \Omega := \R^3 \setminus \overline{B}_1 $, the solution $ u(x) = |x|^{-1} $ satisfies
\[
u \in L^h(\Omega) \; \; \mbox{for} \; \; h \in \left( 3,+\infty \right) \quad \mbox{and} \quad \nabla u \in L^k(\Omega) \; \; \mbox{for} \; \; k \in \left( \frac{3}{2},+\infty \right),
\]
i.e., the summability of the solution and its gradient are different. Hence, the Sobolev space $ W^{1,2}(\R^3) $ is not the natural ambient to look for solutions within, and we will consider spaces which take into account only the summability of the gradient of their elements. Given $ \Omega \subseteq \R^N $, we define (modulo constant functions)
\[
\mathcal{D}^{1,p}(\Omega) := \{ u \in L^1_{\rm loc}(\Omega): \, \nabla u \in L^p(\Omega)\},
\]
endowed with the norm $ \|u\|_{1,p} := \|\nabla u\|_p $ (hereafter $ \|\cdot\|_p $ stands for the standard norm in $ L^p(\Omega) $); then consider its subspace $ \mathcal{D}^{1,p}_0(\Omega) $, called Beppo Levi space, defined as the $ \|\cdot\|_{1,p} $-closure of the set of compactly supported test functions $ C^\infty_{\rm c}(\Omega) $. A Sobolev-type embedding ensures that $ \mathcal{D}^{1,p}_0(\Omega) \hookrightarrow L^{p^*}(\Omega) $. This implies that the functions in $ \mathcal{D}^{1,p}_0(\Omega) $ vanish at infinity, in the sense that for any $ \epsilon > 0 $ the set $ \{x \in \R^N: \, |u(x)| \geq \epsilon \} $ has finite measure. In fact, using also the Chebichev inequality,
\begin{equation}
\label{measdecay}
{\rm meas}(\{x \in \R^N: \, |u(x)| \geq \epsilon \}) \leq \epsilon^{-p^*} \|u\|_{p^*}^{p^*} \leq (c \epsilon^{-1} \|u\|_{1,p})^{p^*} < \infty.
\end{equation}
An open question concerns the possibility to prove, under suitable decay conditions on $ a_i $, that $ u(x), v(x) \to 0 $ as $ |x| \to \infty $; by the way, we will look for solutions to (\ref{problem}) in the product space $ X := \mathcal{D}^{1,p}_0(\R^N) \times \mathcal{D}^{1,q}_0(\R^N) $ (equipped with the norm $ \|(u,v)\|_X := \|u\|_{1,p} + \|v\|_{1,q} $), thus understanding the decay in the measure-theoretic sense described in (\ref{measdecay}).

\bigskip

\noindent
\textbf{Step 1: `Freezing' and `shifting' the right-hand side.}

\medskip

We want to get rid of both singular and convection terms, so we `shift' (i.e., translate by adding a constant) the singular variables of reaction terms by a small quantity, say $ \epsilon > 0 $, and at the same time we `freeze' (i.e., keep fixed) the gradient variables; then we will `unfreeze' the gradient variables by using a fixed point theorem and we will pass to the limit as $ \epsilon \to 0 $, via a priori estimates, to recover a solution of (\ref{problem}). Now we turn into details. \\
Given $ w := (w_1,w_2) \in X \cap C^1_{\rm loc}(\R^N)^2 $ and $ \epsilon > 0 $, we consider the auxiliary problem
\begin{equation}
\label{auxprob1}
\left\{
\begin{array}{lll}
-\Delta_p u = f(x,u+\epsilon,v,\nabla w_1,\nabla w_2) \quad &\mbox{in} \;\; &\R^N, \\
-\Delta_q v = g(x,u,v+\epsilon,\nabla w_1,\nabla w_2) \quad &\mbox{in} \;\; &\R^N, \\
u,v>0 \quad &\mbox{in} \;\; &\R^N,
\end{array}
\right.
\end{equation}
which possesses a unique solution $ (u,v) \in X \cap C^{1,\alpha}_{\rm loc}(\R^N)^2 $, thanks to Minty-Browder's theorem, nonlinear H\"older regularity theory \cite{BCDKS}, and the strong maximum principle. Moreover, setting $ \eta_1 := \max\{\beta_1,\delta_1\} $, $ \eta_2 := \max\{\alpha_2,\gamma_2\} $ and assuming also $ \max\{\|w_i\|_\infty,\|\nabla w_i\|_\infty\} < \infty $, we can deduce the following a priori estimates (see \cite{DM} for $ L^\infty $ estimates on the gradients):
\begin{eqnarray*}
\|\nabla u\|_p^{p-1}& &\leq L_\epsilon (1 + \|\nabla w_1\|_p^{\gamma_1} + \|\nabla w_2\|_q^{\eta_1}), \\
\|\nabla v\|_q^{q-1}& &\leq L_\epsilon (1 + \|\nabla w_1\|_p^{\eta_2} + \|\nabla w_2\|_q^{\delta_2}), \\
\|u\|_\infty& &\leq M_\epsilon(\|\nabla w_1\|_p,\|\nabla w_2\|_q), \\
\|v\|_\infty& &\leq M_\epsilon(\|\nabla w_1\|_p,\|\nabla w_2\|_q), \\
\|\nabla u\|_\infty^{p-1}& &\leq N_\epsilon(\|\nabla w_1\|_p,\|\nabla w_2\|_q,\|w_2\|_\infty) (1 + \|\nabla w_1\|_\infty^{\gamma_1} + \|\nabla w_2\|_\infty^{\delta_1}), \\
\|\nabla v\|_\infty^{q-1}& &\leq N_\epsilon(\|\nabla w_1\|_p,\|\nabla w_2\|_q,\|w_1\|_\infty) (1 + \|\nabla w_1\|_\infty^{\gamma_2} + \|\nabla w_2\|_\infty^{\delta_2}),
\end{eqnarray*}
being $ L_\epsilon > 0 $ a constant and $ M_\epsilon, N_\epsilon $ positive functions which are increasing in each of their arguments.

\bigskip

\noindent
\textbf{Step 2: `Unfreezing' the right-hand side.}

\medskip

In order to `unfreeze' convection terms, that is, to solve the auxiliary problem
\begin{equation}
\label{auxprob2}
\left\{
\begin{array}{lll}
-\Delta_p u = f(x,u+\epsilon,v,\nabla u,\nabla v) \quad &\mbox{in} \;\; &\R^N, \\
-\Delta_q v = g(x,u,v+\epsilon,\nabla u,\nabla v) \quad &\mbox{in} \;\; &\R^N, \\
u,v>0 \quad &\mbox{in} \;\; &\R^N,
\end{array}
\right.
\end{equation}
we consider, for a fixed $ \epsilon > 0 $, the \textit{trapping region}
\begin{equation*}
\begin{split}
\mathcal{R}_\epsilon := \{ &(w_1,w_2) \in X \cap C^1_{\rm loc}(\R^N)^2: \, w_i > 0 \; {\rm in} \; \R^N, \\
&\|\nabla w_1\|_p \leq A_1, \, \|\nabla w_2\|_q \leq A_2, \, \|w_i\|_\infty \leq B, \, \|\nabla w_i\|_\infty \leq C \},
\end{split}
\end{equation*}
with $ A_i,B,C > 0 $ satisfying
\begin{equation}
\label{trapping}
\left\{
\begin{array}{ll}
A_1^{p-1} &\geq L_\epsilon (1 + A_1^{\gamma_1} + A_2^{\eta_1}), \\
A_2^{q-1} &\geq L_\epsilon (1 + A_1^{\eta_2} + A_2^{\delta_2}), \\
B_1 &\geq M_\epsilon(A_1,A_2), \\
B_2 &\geq M_\epsilon(A_1,A_2), \\
C_1^{p-1} &\geq N_\epsilon(A_1,A_2,B_2) (1 + C_1^{\gamma_1} + C_2^{\delta_1}), \\
C_2^{q-1} &\geq N_\epsilon(A_1,A_2,B_1) (1 + C_1^{\gamma_2} + C_2^{\delta_2}).
\end{array}
\right.
\end{equation}
We observe that the algebraic system (\ref{trapping}) is solvable: taking $ 1 < \sigma < \frac{(p-1)(q-1)}{\eta_1 \eta_2} $, it suffices to choose $ A_1 := K^{\frac{1}{\eta_2}} $ and $ A_2 := K^{\frac{\sigma}{q-1}} $, with $ K > 0 $ large enough, to satisfy the first two inequalities, then set $ B_1 := B_2 := M_\epsilon(A_1,A_2) $ to solve the third inequality, and finally pick $ C_1 := H^{\frac{1}{\eta_2}} $, $ C_2 := H^{\frac{\sigma}{q-1}} $, being $ H > 0 $ sufficiently large, to satisfy the last two inequalities; this is essentially due to the $(p,q)$-sublinearity of the right-hand side, although an interplay between the numbers $ A_i,B,C $ occurs. \\
Now we consider the nonlinear operator $ \mathcal{T}_\epsilon(w) := (u,v) $, associating to each $ w \in \mathcal{R}_\epsilon $ the unique solution $ (u,v) $ to (\ref{auxprob1}). According to (\ref{trapping}) we have $ \mathcal{T}(\mathcal{R}_\epsilon) \subseteq \mathcal{R}_\epsilon $; moreover, $ \mathcal{T}_\epsilon $ is continuous and compact in the $ X $-topology, as a consequence of the Sobolev-type embedding, Rellich-Kondrakov's theorem, a diagonal argument, and the $ ({\rm S}_+) $-property of the $p$-Laplacian (in $ \mathcal{D}^{1,p}_0(\R^N) $). Schauder's fixed point theorem then ensures that (\ref{auxprob2}) admits a solution $ (u,v) \in X \cap C^{1,\alpha}_{\rm loc}(\R^N)^2 $. \\
We infer Sobolev and $ L^\infty $ bounds for solutions to (\ref{auxprob2}) by means of (\ref{mixedcond}) and (\ref{weightsum}) respectively, while a Harnack-type result \cite{DaM} yields a local estimate from below. Summarizing, for a large $ C > 0 $ and small $ c_\rho > 0 $ (depending on $ \rho > 0 $) we have
\begin{equation}
\label{properties}
\begin{split}
&\max\{\|\nabla u\|_p,\|\nabla v\|_q,\|u\|_\infty,\|v\|_\infty\} \leq C, \\
&\min\{ {\rm ess \, inf}_{B_\rho} \, u, {\rm ess \, inf}_{B_\rho} \, v \} \geq c_\rho.
\end{split}
\end{equation}

\bigskip

\noindent
\textbf{Step 3: Existence of a distributional solution.}

\medskip

Let us consider a sequence $ \{(u_n,v_n)\} \subseteq X \cap C^1_{\rm loc}(\R^N)^2 $ of solutions to (\ref{auxprob2}) with $ \epsilon = \frac{1}{n} $, and let $ (\varphi_1,\varphi_2) \in C^\infty_{\rm c}(\R^N)^2 $. By reflexivity of $ X $ we get $ (u_n,v_n) \rightharpoonup (u,v) \in X $. We test the equations of (\ref{auxprob2}) with $ (\varphi_1,\varphi_2) $ and then we would like to pass to the limit for $ n \to \infty $ via Lebesgue's dominated convergence theorem. Since $ \varphi_1, \varphi_2 $ are compactly supported and (\ref{properties}) holds for any $ (u_n,v_n) $, solutions are bounded from above and below; so the singular term is not difficult to handle. On the other hand, convection terms are hard to manage, because there would lack compactness on them. Actually, (\ref{properties}) for $ (u_n,v_n) $ ensures that $ f_n(x) := f(\cdot,u_n,v_n,\nabla u_n,\nabla v_n) $ is uniformly bounded in $ L^r(\R^N) $ with $ r > \frac{N}{p} > (p^*)' $, as pointed out in the Introduction; thus, compactness of $ \{\nabla u_n\} $ in $ L^p_{\rm loc}(\R^N) $ is ensured by the following theorem, whose proof is based on the difference quotients method and the Riesz-Fréchet-Kolmogorov compactness criterion.

\textit{
Let $ \Omega \subseteq \R^N $, $ N \geq 2 $, $ p,r \in (1,+\infty) $, and let $ \{u_n\} \in W^{1,p}_{\rm loc}(\Omega) $, $ \{f_n\} \subseteq L^r_{\rm loc}(\Omega) $. Suppose that $ u_n $ is a distributional solution to $ -\Delta_p u_n = f_n $ in $ \Omega $ for all $ n \in \N $. Suppose also that:
\begin{itemize}
	\item $ \{\nabla u_n\} $ is bounded in $ L^p_{\rm loc}(\Omega) $,
	\item $ \{f_n\} $ is bounded in $ L^r_{\rm loc}(\Omega) $,
	\item $ u_n \to u $ in $ L^p_{\rm loc}(\Omega) \cap L^{r'}_{\rm loc}(\Omega) $.
\end{itemize}
Then $ \{\nabla u_n\} $ admits a strongly convergent subsequence in $ L^p_{\rm loc}(\Omega) $.
}

In particular (exploiting also the Rellich-Kondrakov theorem): if $ \{u_n\} $ and $ \{f_n\} $ are uniformly bounded, respectively, in $ W^{1,p}_{\rm loc}(\Omega) $ and $ L^r_{\rm loc}(\Omega) $ with $ r > (p^*)' $, then $ \nabla u_n $ strongly converges (up to subsequences) in $ L^p_{\rm loc}(\Omega) $.

\bigskip

\noindent
\textbf{Step 4: Existence of a weak solution.}

\medskip

Let us consider $ (u,v) \in X $ distributional solution to (\ref{problem}) and take $ \phi \in \mathcal{D}^{1,p}_0(\R^N) $, which can be split into positive and negative parts as $ \phi = \phi^+ - \phi^- $. We perform a localization-regularization procedure, using a set of standard mollifiers $ \{\rho_k\} \subseteq C^\infty_{\rm c}(\R^N) $ and a test function $ \theta \in C^\infty_{\rm c}([0,+\infty)) $ such that
\[
\theta \equiv 1 \;\; \mbox{in} \;\; [0,1], \quad \theta \;\; \mbox{is decreasing in} \;\; (1,2), \quad \theta \equiv 0 \;\; \mbox{in} \;\; [2,+\infty).
\]
We set
\begin{equation*}
\begin{split}
&\theta_n = \theta \left(\frac{|\cdot|}{n}\right) \in C^\infty_{\rm c}(\R^N), \\
&\phi_n = \theta_n \phi^+ \in \mathcal{D}^{1,p}_0(\R^N), \\
&\psi_{k,n} = \rho_k * \phi_n \in C^\infty_{\rm c}(\R^N).
\end{split}
\end{equation*}
We test (\ref{problem}) with $ \psi_{k,n} $ for $ k,n \in \N $; then we pass to the limit in the weak formulation, with respect to both indices, as follows. The properties of mollifiers and (\ref{properties}) allow to pass to the limit in $ k $, while Beppo Levi's monotone convergence theorem and the `good' decay of $ \theta $ (ensuring a high summability of $ \nabla \theta_n $) permit the passage to the limit in $ n $. Repeating the same argument for $ \phi^- $ instead of $ \phi^+ $ proves that $ (u,v) $ is a weak solution to (\ref{problem}).

\begin{small}

\end{small}
\end{document}